\documentclass[11pt]{article}

\usepackage[margin=1.1in]{geometry}
\usepackage{amsmath,amssymb,amsthm}
\usepackage{tikz}
\usepackage{booktabs}
\usepackage{enumitem}
\usepackage{authblk}
\usepackage[hidelinks]{hyperref}

\usetikzlibrary{arrows.meta,decorations.markings,calc}

\newtheorem{theorem}{Theorem}[section]
\newtheorem{lemma}[theorem]{Lemma}
\newtheorem{proposition}[theorem]{Proposition}
\newtheorem{corollary}[theorem]{Corollary}

\theoremstyle{definition}
\newtheorem{definition}[theorem]{Definition}
\newtheorem{example}[theorem]{Example}
\theoremstyle{remark}
\newtheorem{remark}[theorem]{Remark}
\newtheorem*{question}{Question}
\theoremstyle{plain}
\newtheorem{thmx}{Theorem}

\newtheorem{conjx}[thmx]{Conjecture}

\newcommand{\E}{\mathcal{E}}
\newcommand{\T}{\mathcal{T}}
\newcommand{\C}{\mathcal{C}}
\newcommand{\Cx}{\mathbb{C}}
\newcommand{\Z}{\mathbb{Z}}
\newcommand{\fc}{\mathrm{f}}       
\newcommand{\Up}{U}
\newcommand{\Dn}{D}

\tikzset{
  eis/.style={x={(1cm,0cm)}, y={(0.5cm,0.8660254cm)}},
  blk/.style={fill=black!55},
  wht/.style={fill=white},
  grid/.style={draw=black!45, line width=0.3pt},
  bd/.style={draw=black, line width=0.9pt},
}

\newcommand{\DrawDelta}[3]{%
  \pgfmathtruncatemacro{\dn}{#1}%
  \pgfmathtruncatemacro{\dnm}{#1-1}%
  \foreach \i in {0,...,\dnm}{%
    \pgfmathtruncatemacro{\jmax}{\dn-\i-1}%
    \foreach \j in {0,...,\jmax}{%
      \path[#2] (\i,\j) -- ++(1,0) -- ++(-1,1) -- cycle;%
    }%
    \ifnum\jmax>0
      \pgfmathtruncatemacro{\jmaxd}{\jmax-1}%
      \foreach \j in {0,...,\jmaxd}{%
        \path[#3] (\i+1,\j) -- ++(-1,1) -- ++(1,0) -- cycle;%
      }%
    \fi
  }%
  \ifnum\dn>1
    \foreach \k in {1,...,\dnm}{%
      \pgfmathtruncatemacro{\r}{\dn-\k}%
      \draw[grid] (\k,0) -- (\k,\r);%
      \draw[grid] (0,\k) -- (\r,\k);%
      \draw[grid] (\k,0) -- (0,\k);%
    }%
  \fi
  \draw[bd] (0,0) -- (\dn,0) -- (0,\dn) -- cycle;%
}

\begin{document}

\title{Good colorings of doubled triangles\\ and a conjecture of Schwartz}

\author[1]{Lakshitha Fernando}
\author[2,*]{Eranda Dhananjaya}
\author[2]{Rajitha Ranasinghe}
\author[3]{Dinusha Priyadarshani}

\affil[1]{School of Data, Mathematical, and Statistical Sciences, University of
Central Florida, FL, USA}
\affil[2]{Department of Mathematics, Faculty of Science, University of
Peradeniya, Peradeniya, Sri Lanka}
\affil[3]{Department of Basic Sciences, Faculty of Health Sciences, The Open
University of Sri Lanka, Nawala, Sri Lanka}
\affil[*]{Corresponding author. Email: \texttt{dhananjayae@sci.pdn.ac.lk}}

\date{}

\maketitle

\begin{abstract}
A good coloring of a sphere triangulation is a black-and-white coloring
of its faces such that, at every vertex, the numbers of black and white
faces differ by a multiple of~$3$. Good colorings correspond to proper
vertex $4$-colorings. For a good coloring with fold count $\mathrm{f}$
(the number of black-white edges) on a triangulation with $F$ faces,
R.~E.~Schwartz introduced the Eisenstein isoperimetric ratio
$\eta=\mathrm{f}^2/F$. For triangulations with degree sequence
$6,\dots,6,2,2,2$, he proved that $\eta\ge3$ and conjectured that every
such triangulation has a good coloring with $\eta$ below a universal
constant. These triangulations are indexed by the nonzero Eisenstein
integers~$\beta$. We study the family $\T(n)$ given by $\beta=n$, a
positive integer, which consists of the doubled equilateral triangles of
side~$n$. We prove that every good coloring of $\T(n)$ has at least $3n$
folds, so $\eta\ge9/2$, and that this bound is attained for every~$n$.
Thus the minimum of $\eta$ on this family is exactly $9/2$, and
Schwartz's conjecture holds for it. We also prove that subdividing a
good coloring of $\T(\beta)$ gives a good coloring of $\T(n\beta)$ with
the same value of $\eta$. Consequently, Schwartz's conjecture holds in
general if and only if it holds for primitive Eisenstein integers.
\end{abstract}



\section{Introduction}\label{sec:intro}

\subsection{Background}

The Four Color Theorem~\cite{AH1,AHK2,Gonthier,RSST} states that the
vertices of every triangulation of the sphere can be properly colored
with four colors; see~\cite{Wilson} for its history. Label the vertices
of a regular tetrahedron by the four colors. A proper vertex
$4$-coloring of a triangulation $\Sigma$ then induces a simplicial map
from $\Sigma$ to the tetrahedron. Color a face of $\Sigma$ black if this
map preserves orientation on it, and white otherwise. Since three faces
of the tetrahedron meet at each vertex, at every vertex of $\Sigma$ the
number of black faces is congruent to the number of white faces
modulo~$3$. Such a face coloring is called a \emph{good coloring}.
Conversely, every good coloring of a sphere triangulation comes from a
proper vertex $4$-coloring~\cite[\S1.1]{Sch}. This point of view is
close to the geometric coloring theory of Fisk~\cite{Fisk}.

Schwartz~\cite{Sch} used it to study the geometry of good colorings on
the simplest triangulations of non-negative combinatorial curvature.
These are the triangulations with degree sequence $6,\dots,6,2,2,2$, that
is, three vertices of degree~$2$ and all other vertices of degree~$6$.
(We allow the case with no vertex of degree~$6$.) Following
Schwartz~\cite{Sch} and Thurston~\cite{Thurston}, they are indexed by
the nonzero Eisenstein integers~$\beta$. We write $\T(\beta)$ for the
triangulation that corresponds to~$\beta$. For a good coloring $\C$, the
\emph{fold count} $\fc(\C)$ is the number of edges that separate a black
face from a white face. Schwartz defined the \emph{Eisenstein
isoperimetric ratio}
\[
  \eta(\C) = \frac{\fc(\C)^2}{F},
\]
where $F$ is the number of faces. It is small when the black and white
regions are large and ``round''. He proved the following lower bound.

\begin{thmx}[{\cite[Theorem~1.2]{Sch}}]\label{thm:sch-lower}
Let $\C$ be a good coloring of a triangulation with degree sequence
$6,\dots,6,2,2,2$. Then $\eta(\C)\ge 3$.
\end{thmx}

He also made two conjectures. Let $\alpha=(1+i\sqrt3)/2$, so that the
Eisenstein integers are the numbers $a+b\alpha$ with $a,b\in\Z$. A
triangulation $\T(\beta)$ with $\beta=a+b\alpha$ is called
\emph{primitive} if $\gcd(a,b)=1$. Let
$\phi=(1+\sqrt5)/2$ be the golden ratio.

\begin{conjx}[{\cite[Conjecture~1.3]{Sch}}]\label{conj:phi}
Let $\{\mathcal G_m\}$ be an infinite sequence of distinct good colorings of
primitive triangulations with degree sequence $6,\dots,6,2,2,2$. Then
$\liminf_{m\to\infty}\eta(\mathcal G_m)\ge\phi^6$.
\end{conjx}

\begin{conjx}[{\cite[Conjecture~1.4]{Sch}}]\label{conj:omega}
There is a constant $\Omega$ with the following property. For every
triangulation $\Sigma$ with degree sequence $6,\dots,6,2,2,2$ there is a
good coloring $\C$ of $\Sigma$ with $\eta(\C)<\Omega$.
\end{conjx}

To our knowledge, Conjecture~\ref{conj:omega} is open. Schwartz's
continued fraction colorings have bounded $\eta$ along the continued
fraction approximants of a quadratic irrational~\cite[Theorem~1.1]{Sch},
but for other parameters, such as $1/13$, their fold count is linear in
the number of faces~\cite[\S2.5]{Sch}. Apart
from the general bound $\eta\ge3$, we are not aware of an infinite family
of these triangulations on which the minimum of $\eta$ over all good
colorings is known.

\subsection{Our contributions}

We study the family $\T(n)$, where $\beta=n$ is a positive integer.
Geometrically, $\T(n)$ is a doubled equilateral triangle of side~$n$ with
$2n^2$ faces. Our contributions are the following.

\smallskip
\noindent\textbf{(1) The exact minimum of $\eta$ on $\T(n)$.}

\begin{theorem}\label{thm:main}
Let $n\ge1$. Every good coloring $\C$ of $\T(n)$ satisfies
$\fc(\C)\ge 3n$, and hence
\[
  \eta(\C)\ \ge\ \frac{9}{2}.
\]
Moreover, $\T(n)$ has a good coloring $\mathcal B_n$ with
$\fc(\mathcal B_n)=3n$ and $\eta(\mathcal B_n)=9/2$.
\end{theorem}

In $\mathcal B_n$, one of the two triangles is black and the other is
white (see Figure~\ref{fig:AB}). The substance of
Theorem~\ref{thm:main} is the lower bound, which holds for \emph{every}
good coloring of $\T(n)$. It improves the bound of
Theorem~\ref{thm:sch-lower} from $3$ to $9/2$ on this family, and it is
sharp for every~$n$. The proof is different from Schwartz's argument for
Theorem~\ref{thm:sch-lower}: we cut $\T(n)$ into $3n$ closed bands of
triangles, show that each band contains at least two folds, and count.
For comparison, the checkerboard coloring $\mathcal A_n$ of $\T(n)$ has
every edge as a fold, so $\eta(\mathcal A_n)=9n^2/2$
(Proposition~\ref{prop:A}), which is unbounded.

\smallskip
\noindent\textbf{(2) Invariance of $\eta$ under subdivision.}

\begin{theorem}\label{thm:scaling}
Let $\beta$ be a nonzero Eisenstein integer and let $n\ge1$ be an integer. If $\C$ is
a good coloring of $\T(\beta)$, then there is a good coloring $\C^{(n)}$
of $\T(n\beta)$ with
\[
  \fc\bigl(\C^{(n)}\bigr)=n\,\fc(\C)
  \quad\text{and}\quad
  \eta\bigl(\C^{(n)}\bigr)=\eta(\C).
\]
\end{theorem}

The coloring $\mathcal B_n$ is $\C^{(n)}$ for the good coloring
$\C$ of $\T(1)$, which is unique up to swapping colors.

\smallskip
\noindent\textbf{(3) Consequences for Schwartz's conjectures.}

\begin{corollary}\label{cor:reduction}
\leavevmode
\begin{enumerate}[label=\textup{(\roman*)}]
\item Conjecture~\ref{conj:omega} holds for the family
  $\{\T(n):n\ge1\}$ with any constant $\Omega>9/2$.
\item Conjecture~\ref{conj:omega} holds if and only if it holds, with the same
  constant $\Omega$, for all primitive triangulations $\T(\beta)$.
\end{enumerate}
\end{corollary}

\begin{remark}\label{rem:primitive}
For $n\ge2$ the triangulation $\T(n)$ is not primitive. Since
$9/2<\phi^6\approx17.94$, Proposition~\ref{prop:B} shows that the word
``primitive'' cannot be removed from Conjecture~\ref{conj:phi}. Indeed,
the colorings $\mathcal B_2,\mathcal B_3,\dots$ are distinct and all have
$\eta=9/2$. More generally, for any $\beta$ and any good coloring $\C$
of $\T(\beta)$, Theorem~\ref{thm:scaling} gives infinitely many distinct
colorings $\C^{(n)}$ of the triangulations $\T(n\beta)$, all with the
same value of $\eta$.
\end{remark}

Part~(ii) of Corollary~\ref{cor:reduction} shows that
Conjecture~\ref{conj:omega} is really a question about primitive
triangulations, and Remark~\ref{rem:primitive} shows that the
restriction to primitive triangulations in Conjecture~\ref{conj:phi} is
necessary.

\subsection{Organization}

In Section~\ref{sec:prelim} we fix the definitions and describe $\T(n)$
as a doubled triangle. In Section~\ref{sec:examples} we study the
colorings $\mathcal A_n$ and $\mathcal B_n$ and give some examples.
Section~\ref{sec:lower} proves the lower bound in
Theorem~\ref{thm:main}. Section~\ref{sec:scaling} proves
Theorem~\ref{thm:scaling} and Corollary~\ref{cor:reduction}. We end
with some remarks and questions.

\section{Preliminaries}\label{sec:prelim}

\subsection{Triangulations and good colorings}

As in~\cite{Sch,Thurston}, a \emph{triangulation} of the sphere is a
decomposition of the sphere into finitely many triangles, glued edge to
edge. We allow two triangles to share more than one edge. This is
needed because the two triangles at a vertex of degree~$2$ always share
both of their edges at that vertex. The \emph{degree} of a vertex is the
number of triangle corners at that vertex. We write $V$, $E$ and $F$
for the numbers of vertices, edges and faces. Since every face has
three edges and every edge lies on two faces, $2E=3F$. Euler's formula
gives $V-E+F=2$.

\begin{definition}\label{def:good}
A \emph{good coloring} of a triangulation is a coloring of its faces in
black and white such that, at every vertex $v$,
\[
  b(v)-w(v)\equiv 0 \pmod 3,
\]
where $b(v)$ and $w(v)$ are the numbers of black and white face corners
at~$v$. The \emph{fold count} $\fc(\C)$ of a face coloring $\C$ is the
number of edges whose two sides have different colors. Such an edge is
called a \emph{fold}. The \emph{Eisenstein isoperimetric ratio} is
$\eta(\C)=\fc(\C)^2/F$.
\end{definition}

We will often use the following simple reformulation. If $v$ has degree
$d$, then $w(v)=d-b(v)$, so the condition reads $2b(v)\equiv d
\pmod 3$. Multiplying by $2$ gives $b(v)\equiv 2d\pmod3$.

\begin{lemma}\label{lem:local}
Let $\C$ be a face coloring of a triangulation.
\begin{enumerate}[label=\textup{(\alph*)}]
\item At a vertex of degree~$6$, the condition of
  Definition~\ref{def:good} holds if and only if $b(v)\in\{0,3,6\}$.
\item At a vertex of degree~$2$, it holds if and only if $b(v)=1$. In
  other words, the two faces at $v$ have different colors.
\end{enumerate}
\end{lemma}

\begin{proof}
If $d=6$ then $b(v)\equiv 12\equiv 0\pmod 3$, and $0\le b(v)\le 6$. If
$d=2$ then $b(v)\equiv 4\equiv1\pmod3$, and $0\le b(v)\le2$.
\end{proof}

Swapping black and white turns a good coloring into a good coloring
with the same fold count. We will use this symmetry without comment.

\subsection{The triangulations \texorpdfstring{$\T(\beta)$}{T(beta)}}

Let $\alpha=(1+i\sqrt3)/2$ and let $\E=\{a+b\alpha:a,b\in\Z\}$ be the
ring of Eisenstein integers. We follow the notation of~\cite{Sch}. The
points of $\E$ are the vertices of the regular triangulation $\T$ of
the plane $\Cx$ by unit equilateral triangles. For a nonzero
$\beta\in\E$, let $G_\beta$ be the group generated by the rotations of
order~$3$ about the points of $\beta\E$. Then
\[
  \T(\beta)=\T/G_\beta
\]
is a triangulation of the sphere. It has three vertices of degree~$2$,
which are the images of the points of $\beta\E$. All other vertices
have degree~$6$. A fundamental domain for $G_\beta$ is the rhombus
$R_\beta$ with vertices $0$, $\beta$, $\beta+\beta\alpha$ and
$\beta\alpha$~\cite[\S1.2]{Sch}. The rotation of order~$3$ about
$\beta$ maps the side $[\beta,0]$ onto the side $[\beta,\beta+\beta\alpha]$.
The rotation of order~$3$ about $\beta\alpha$ maps the side
$[\beta\alpha,0]$ onto the side $[\beta\alpha,\beta+\beta\alpha]$. The
sphere $\T(\beta)$ is obtained from $R_\beta$ by these two gluings.
Every triangulation of the sphere with degree sequence
$6,\dots,6,2,2,2$ is isomorphic to $\T(\beta)$ for some
nonzero $\beta$~\cite{Sch,Thurston}. The number of faces is
$F=2|\beta|^2=2(a^2+ab+b^2)$ for $\beta=a+b\alpha$.

\subsection{The doubled triangle}\label{subsec:doubled}

We now describe $\T(n)$ in a way that is easy to work with. We use the
coordinates $(x,y)$ for the point $x+y\alpha$ of the plane. Let
\[
  \Delta_n=\{(x,y): x\ge0,\ y\ge0,\ x+y\le n\}.
\]
This is an equilateral triangle of side $n$ with vertices $(0,0)$,
$(n,0)$ and $(0,n)$. The lattice $\T$ cuts it into $n^2$ unit
triangles. For integers $i,j\ge0$ we write
\begin{align*}
  \Up(i,j) &= \text{the triangle with vertices } (i,j),\ (i+1,j),\ (i,j+1),
     && i+j\le n-1,\\
  \Dn(i,j) &= \text{the triangle with vertices } (i+1,j),\ (i,j+1),\ (i+1,j+1),
     && i+j\le n-2.
\end{align*}
We call them \emph{up} and \emph{down} triangles. There are $n(n+1)/2$
up triangles and $n(n-1)/2$ down triangles.

Take two copies $\Delta_n^+$ and $\Delta_n^-$ of $\Delta_n$. Glue them
along their boundaries by the identity map. The result is a sphere
$S_n$, triangulated by the $2n^2$ unit triangles $\Up^\pm(i,j)$ and
$\Dn^\pm(i,j)$. We call $\Delta_n^+$ the \emph{front} and $\Delta_n^-$
the \emph{back}. A point on the boundary of $\Delta_n$ is a single point
of $S_n$. A point in the interior of $\Delta_n$ gives two points of
$S_n$, one on each copy. The same holds for edges. We write
$P_0=(0,0)$, $P_1=(n,0)$ and $P_2=(0,n)$ for the three corners.

\begin{lemma}\label{lem:Sn}
The triangulation $S_n$ is isomorphic to $\T(n)$.
\end{lemma}

\begin{proof}
For $\beta=n$ the rhombus $R_n$ has vertices $(0,0)$, $(n,0)$, $(n,n)$
and $(0,n)$. It is the union of $\Delta_n$ and the triangle $\Delta'_n$
with vertices $(n,0)$, $(n,n)$, $(0,n)$. The two triangles share the
side $[(n,0),(0,n)]$. Let $\psi$ be the reflection of the plane in the
line through $(n,0)$ and $(0,n)$. It maps the lattice $\T$ to itself,
maps $\Delta'_n$ onto $\Delta_n$, fixes the common side, and sends
$(n,n)$ to $(0,0)$.

Define a map $h\colon R_n\to S_n$ as follows. On $\Delta_n$ it is the
identity onto the front $\Delta_n^+$. On $\Delta'_n$ it is $\psi$,
followed by the identity onto the back $\Delta_n^-$. The two
definitions agree on the common side, because $\psi$ fixes it and the
side lies on $\partial\Delta_n$. The map $h$ sends unit triangles to
unit triangles.

Let ${\sim}$ be the equivalence relation on $R_n$ generated by the two
gluings described above, so that $\T(n)=R_n/{\sim}$. We check that $h$
respects these gluings. The rotation about
$(n,0)$ glues the point at distance $t$ from $(n,0)$ on $[(n,0),(0,0)]$
to the point at distance $t$ from $(n,0)$ on $[(n,0),(n,n)]$. The map
$\psi$ sends the second point to the point at distance $t$ from $(n,0)$
on $[(n,0),(0,0)]$. In $S_n$ the front and back copies of this point
are the same point. So $h$ maps the two glued points to the same point.
The same argument works for the gluing at $(0,n)$. Hence $h$ induces a
bijection from $R_n/{\sim}$ onto $S_n$ that maps triangles to
triangles. This is an isomorphism of triangulations.
\end{proof}

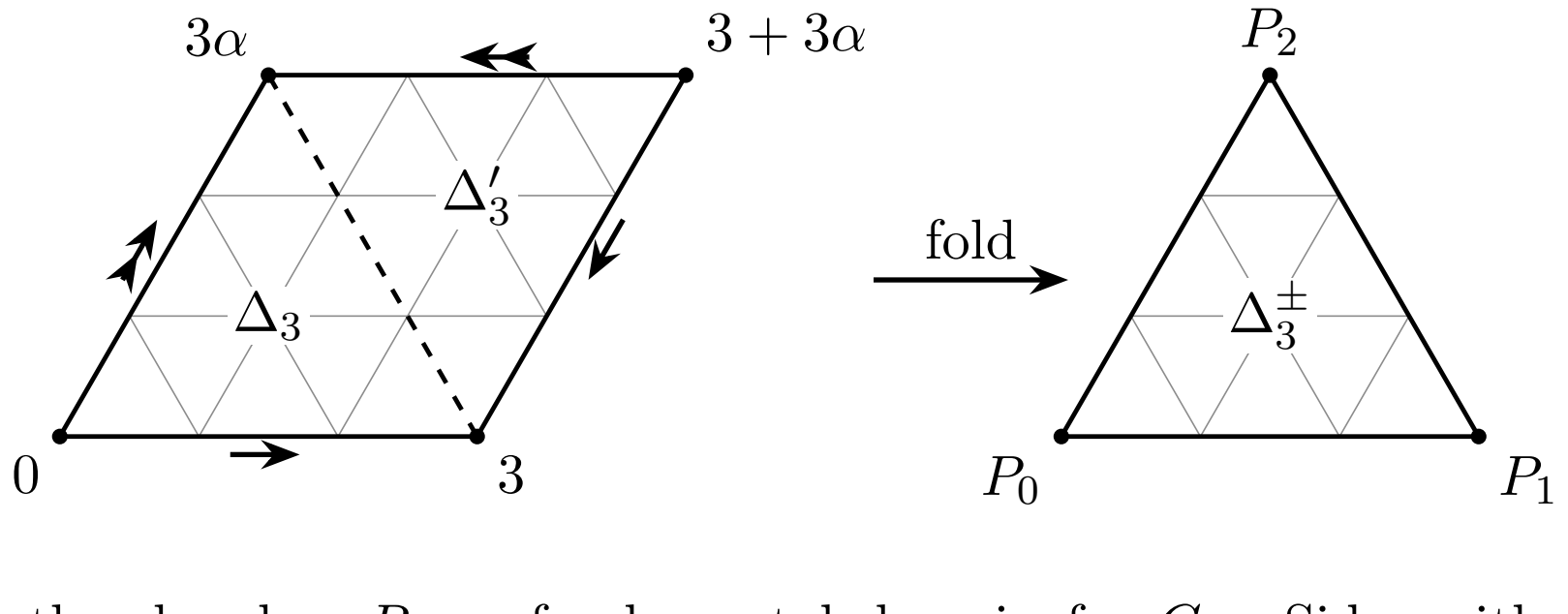
\begin{figure}[ht]
\centering
\begin{tikzpicture}[eis, scale=0.95]
  \begin{scope}
    \foreach \k in {1,2}{
      \draw[grid] (\k,0) -- (\k,3);
      \draw[grid] (0,\k) -- (3,\k);
      \draw[grid] (\k,0) -- (0,\k);
      \draw[grid] (3,\k) -- (\k,3);
    }
    \draw[bd] (0,0) -- (3,0);
    \draw[bd] (0,0) -- (0,3);
    \draw[bd] (3,0) -- (3,3) -- (0,3);
    \draw[bd, dashed] (3,0) -- (0,3);
    \draw[line width=1pt,-{Stealth}] (1.3,-0.15) -- (1.8,-0.15);
    \draw[line width=1pt,-{Stealth}] (3.15,1.8) -- (3.15,1.3);
    \draw[line width=1pt,-{Stealth[]Stealth[]}] (-0.2,1.3) -- (-0.2,1.8);
    \draw[line width=1pt,-{Stealth[]Stealth[]}] (1.8,3.15) -- (1.3,3.15);
    \fill (0,0) circle (1.6pt) node[below left] {$0$};
    \fill (3,0) circle (1.6pt) node[below right] {$3$};
    \fill (0,3) circle (1.6pt) node[above left] {$3\alpha$};
    \fill (3,3) circle (1.6pt) node[above right] {$3+3\alpha$};
    \node[fill=white,inner sep=1pt] at (1,1) {$\Delta_3$};
    \node[fill=white,inner sep=1pt] at (2,2) {$\Delta'_3$};
  \end{scope}
  \draw[-{Stealth}, line width=1pt] (5.2,1.3) -- (6.6,1.3)
     node[midway, above] {fold};
  \begin{scope}[shift={(7.2,0)}]
    \DrawDelta{3}{wht}{wht}
    \fill (0,0) circle (1.6pt) node[below left] {$P_0$};
    \fill (3,0) circle (1.6pt) node[below right] {$P_1$};
    \fill (0,3) circle (1.6pt) node[above] {$P_2$};
    \node[fill=white,inner sep=1pt] at (1,1) {$\Delta_3^{\pm}$};
  \end{scope}
\end{tikzpicture}
\caption{Left: the rhombus $R_3$, a fundamental domain for $G_3$. Sides
with the same arrow are glued by a rotation of order~$3$. Right: folding
$\Delta'_3$ onto $\Delta_3$ along the dashed side gives the doubled
triangle $S_3=\T(3)$. The three corners $P_0,P_1,P_2$ are the
vertices of degree~$2$.}
\label{fig:rhombus}
\end{figure}

From now on we identify $\T(n)$ with $S_n$. Figure~\ref{fig:rhombus}
shows the rhombus and the folding for $n=3$.

\begin{lemma}\label{lem:counts}
The triangulation $\T(n)$ has $F=2n^2$ faces, $E=3n^2$ edges and
$V=n^2+2$ vertices. Exactly three vertices have degree~$2$, namely the
corners $P_0$, $P_1$ and $P_2$. The other $n^2-1$
vertices have degree~$6$.
\end{lemma}

\begin{proof}
We count in $S_n$. There are $2n^2$ faces, so $E=3F/2=3n^2$, and
$V=2-F+E=n^2+2$. A corner of $\Delta_n$ lies on one unit triangle of
each copy, so it has degree~$2$. A point of $\partial\Delta_n$ that is
not a corner lies on three unit triangles of each copy, so it has
degree~$6$. An interior lattice point of $\Delta_n$ lies on six unit
triangles of its own copy, so it also has degree~$6$. Hence there are
$V-3=n^2-1$ vertices of degree~$6$.

We can also count them directly. Each side of $\Delta_n$ has $n-1$
lattice points that are not corners. The interior of $\Delta_n$ has
$(n-1)(n-2)/2$ lattice points, and each one appears twice in $S_n$. So
the number of vertices of degree~$6$ is
\[
  3(n-1)+2\cdot\frac{(n-1)(n-2)}{2}=(n-1)(n+1)=n^2-1. \qedhere
\]
\end{proof}

\section{Two good colorings and some examples}\label{sec:examples}

\begin{figure}[ht]
\centering
\begin{tikzpicture}[eis, scale=0.8]
  \begin{scope}
    \DrawDelta{3}{blk}{blk}
    \node at (1.5,-0.55) {front $\Delta_3^+$};
  \end{scope}
  \begin{scope}[shift={(4,0)}]
    \DrawDelta{3}{wht}{wht}
    \node at (1.5,-0.55) {back $\Delta_3^-$};
  \end{scope}
  \node at (-1.3,1.0) {$\mathcal B_3$};
  \begin{scope}[shift={(2.2,-4.4)}]
    \DrawDelta{3}{blk}{wht}
    \node at (1.5,-0.55) {front $\Delta_3^+$};
  \end{scope}
  \begin{scope}[shift={(6.2,-4.4)}]
    \DrawDelta{3}{wht}{blk}
    \node at (1.5,-0.55) {back $\Delta_3^-$};
  \end{scope}
  \node at (0.9,-3.4) {$\mathcal A_3$};
\end{tikzpicture}
\caption{Top: the two-sided coloring $\mathcal B_3$ of $\T(3)$. Its
folds are the $9$ edges on the common boundary, and $\eta=81/18=9/2$.
Bottom: the checkerboard coloring $\mathcal A_3$. All $27$ edges are
folds, and $\eta=729/18=81/2$. Black faces are drawn in gray so that
the grid stays visible.}
\label{fig:AB}
\end{figure}
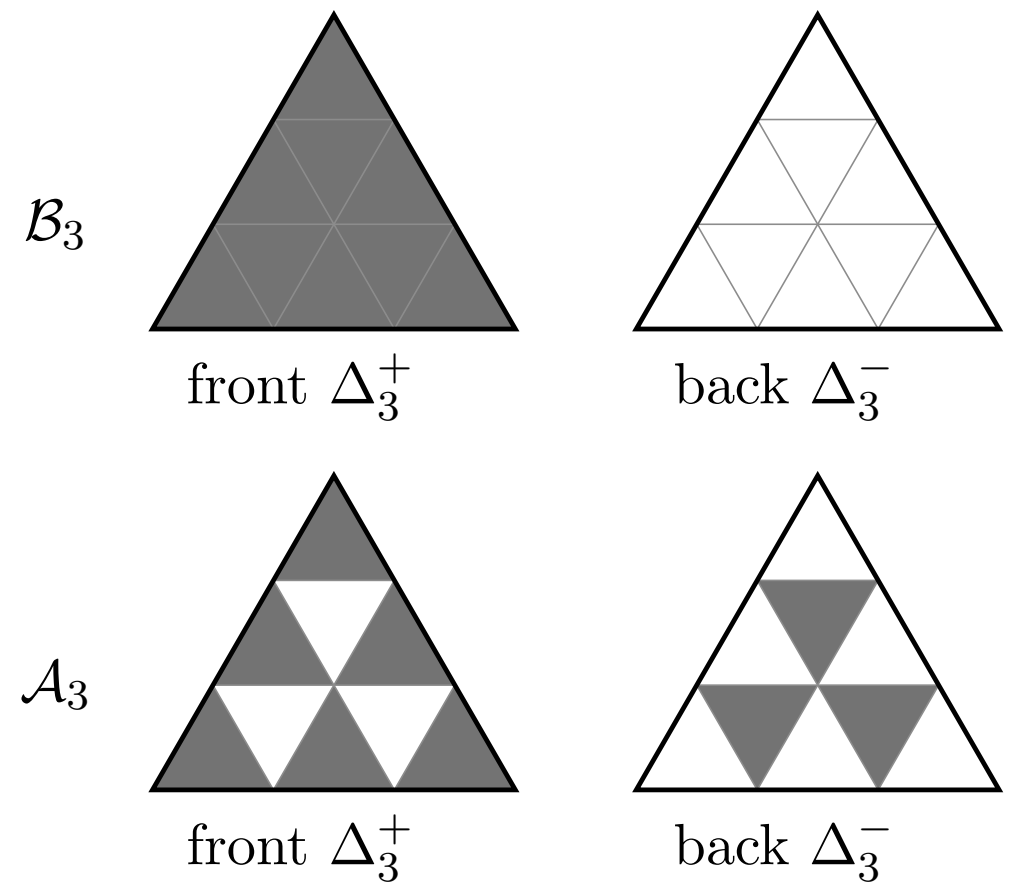

We now define two good colorings of $\T(n)=S_n$. Both are shown in
Figure~\ref{fig:AB}.

\begin{definition}
We define two face colorings of $S_n$.
\begin{enumerate}[label=\textup{(\roman*)}]
\item The \emph{checkerboard coloring} $\mathcal A_n$ colors
  $\Up^+(i,j)$ and $\Dn^-(i,j)$ black, and $\Dn^+(i,j)$ and
  $\Up^-(i,j)$ white.
\item The \emph{two-sided coloring} $\mathcal B_n$ colors every face of
  the front $\Delta_n^+$ black and every face of the back $\Delta_n^-$
  white.
\end{enumerate}
\end{definition}

\begin{proposition}\label{prop:A}
The coloring $\mathcal A_n$ is good. Every edge of $\T(n)$ is a fold of
$\mathcal A_n$. Hence $\fc(\mathcal A_n)=3n^2$ and
$\eta(\mathcal A_n)=9n^2/2$.
\end{proposition}

\begin{proof}
Let $e$ be an edge. If $e$ lies inside one copy of $\Delta_n$, then
its two sides are an up triangle and a down triangle of the same copy.
These have different colors. If $e$ lies on $\partial\Delta_n$, then
its two sides are the front and back copies of the same up triangle.
(Every unit triangle of $\Delta_n$ with an edge on $\partial\Delta_n$
is an up triangle.) These also have different colors. So every edge is
a fold. Going around a vertex, the colors therefore alternate, so
$b(v)=w(v)$ at every vertex, and $\mathcal A_n$ is good. Finally,
$\fc(\mathcal A_n)=E=3n^2$ and $\eta=9n^4/(2n^2)=9n^2/2$.
\end{proof}

The checkerboard coloring is the ``alternating coloring'' mentioned
in~\cite[\S1.2]{Sch}. It has the largest possible fold count, so it is the
worst coloring from the point of view of $\eta$.

\begin{proposition}\label{prop:B}
The coloring $\mathcal B_n$ is good, $\fc(\mathcal B_n)=3n$, and
$\eta(\mathcal B_n)=9/2$.
\end{proposition}

\begin{proof}
We use Lemma~\ref{lem:local}. A corner $P_k$ lies on one front face and
one back face, so $b(P_k)=1$. A vertex on $\partial\Delta_n$ that is not
a corner lies on three front faces and three back faces, so $b(v)=3$.
A vertex in the interior of the front has $b(v)=6$. A vertex in the
interior of the back has $b(v)=0$. So $\mathcal B_n$ is good.

An edge inside $\Delta_n^\pm$ has two faces of the same color. An edge
on $\partial\Delta_n$ has one front face and one back face. So the folds
are exactly the $3n$ edges on $\partial\Delta_n$. Hence
$\eta(\mathcal B_n)=(3n)^2/(2n^2)=9/2$.
\end{proof}

\begin{example}[The case $n=1$]
The triangulation $\T(1)$ has two faces and three vertices, each of
degree~$2$. By Lemma~\ref{lem:local}(b), the two faces must have
different colors. So $\T(1)$ has exactly one good coloring up to
swapping colors, and $\mathcal A_1=\mathcal B_1$. It has $\fc=3$ and
$\eta=9/2$.
\end{example}

\begin{figure}[ht]
\centering
\begin{tikzpicture}[eis, scale=1.35]
  \begin{scope}
    \DrawDelta{2}{blk}{blk}
    \foreach \p/\c in {{(0,0)}/4,{(2,0)}/4,{(0,2)}/4,{(1,0)}/1,{(1,1)}/2,{(0,1)}/3}{
      \fill[white] \p circle (5pt); \draw \p circle (5pt);
      \node[font=\small] at \p {\c};
    }
    \node at (1,-0.5) {front};
  \end{scope}
  \begin{scope}[shift={(3.2,0)}]
    \DrawDelta{2}{wht}{wht}
    \foreach \p/\c in {{(0,0)}/4,{(2,0)}/4,{(0,2)}/4,{(1,0)}/1,{(1,1)}/2,{(0,1)}/3}{
      \fill[white] \p circle (5pt); \draw \p circle (5pt);
      \node[font=\small] at \p {\c};
    }
    \node at (1,-0.5) {back};
  \end{scope}
\end{tikzpicture}
\caption{A proper vertex $4$-coloring of $\T(2)$. The front is a net of
the tetrahedron. The induced good coloring is $\mathcal B_2$.}
\label{fig:tetra}
\end{figure}
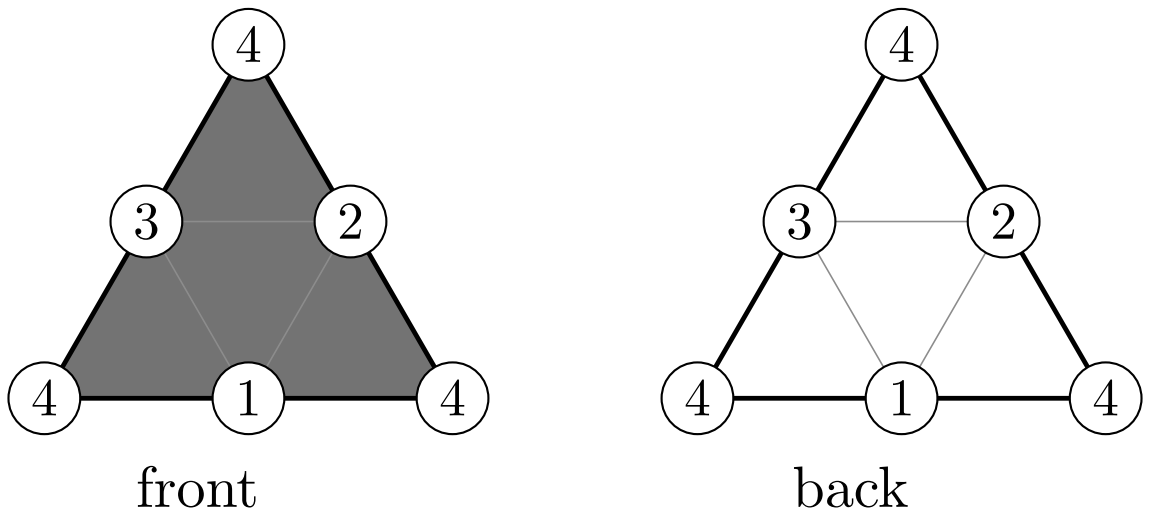

\begin{example}[The case $n=2$ and the tetrahedron]\label{ex:n2}
For $n=2$ the front $\Delta_2^+$ is cut into four unit triangles. This
is exactly the usual net of a regular tetrahedron. Color the three
corners with color $4$ and the three midpoints of the sides with colors
$1$, $2$, $3$, as in Figure~\ref{fig:tetra}. This is a proper vertex
$4$-coloring of $\T(2)$. The induced map to the tetrahedron folds the
front onto the surface of the tetrahedron, preserving orientation on
every front face. The back faces carry the same labels, but they have
the opposite orientation in $S_2$. So the induced good coloring is
$\mathcal B_2$ (up to swapping colors). Here $\fc=6$, $F=8$ and
$\eta=36/8=9/2$.
\end{example}

Table~\ref{tab:values} compares the two colorings for small $n$.

\begin{table}[!ht]
\centering
\begin{tabular}{@{}cccccc@{}}
\toprule
$n$ & $F=2n^2$ & $\fc(\mathcal A_n)$ & $\eta(\mathcal A_n)$
    & $\fc(\mathcal B_n)$ & $\eta(\mathcal B_n)$\\
\midrule
1 & 2  & 3  & 4.5   & 3  & 4.5\\
2 & 8  & 12 & 18    & 6  & 4.5\\
3 & 18 & 27 & 40.5  & 9  & 4.5\\
4 & 32 & 48 & 72    & 12 & 4.5\\
5 & 50 & 75 & 112.5 & 15 & 4.5\\
\bottomrule
\end{tabular}
\medskip
\caption{Fold counts and isoperimetric ratios of $\mathcal A_n$ and
$\mathcal B_n$.}
\label{tab:values}
\end{table}

\section{The lower bound}\label{sec:lower}

In this section we prove that every good coloring of $\T(n)$ has at
least $3n$ folds. The idea is to cut $S_n$ into closed bands of
triangles, which we call rings. We show that each ring contains at least
two folds. Then we count.

\subsection{Rings}

Put $z=n-x-y$. On $\Delta_n$ the three functions $x$, $y$, $z$ are
non-negative and their sum is~$n$. Each side of $\Delta_n$ is the zero
set of one of them. Along every edge of the lattice, exactly one of
$x,y,z$ is constant, and the other two change by exactly~$1$.

Fix $k\in\{0,1,\dots,n-1\}$. The strip $\{k\le x\le k+1\}\cap\Delta_n$
contains the $2(n-k)-1$ unit triangles
\[
  \Up(k,0),\ \Dn(k,0),\ \Up(k,1),\ \Dn(k,1),\ \dots,\
  \Dn(k,n-k-2),\ \Up(k,n-k-1),
\]
listed in this order. Two consecutive triangles in the list share an
edge. The first triangle $\Up(k,0)$ has an edge on the side $y=0$, and
the last triangle $\Up(k,n-k-1)$ has an edge on the side $z=0$.

\begin{definition}
The \emph{$x$-ring} $X_k$ is the cyclic sequence of faces of $S_n$
\[
  \Up^+(k,0),\ \Dn^+(k,0),\ \dots,\ \Up^+(k,n-k-1),\
  \Up^-(k,n-k-1),\ \dots,\ \Dn^-(k,0),\ \Up^-(k,0),
\]
which then closes up at $\Up^+(k,0)$. A \emph{crossing edge} of $X_k$
is an edge of $S_n$ that lies in the strip $\{k\le x\le k+1\}$ (on the
front or on the back) and along which $x$ is not constant.
\end{definition}

\begin{figure}[ht]
\centering
\begin{tikzpicture}[eis, scale=0.95]
  \DrawDelta{5}{wht}{wht}
  \foreach \j in {0,...,3}{
    \path[fill=black!35] (1,\j) -- ++(1,0) -- ++(-1,1) -- cycle;
  }
  \foreach \j in {0,...,2}{
    \path[fill=black!35] (2,\j) -- ++(-1,1) -- ++(1,0) -- cycle;
  }
  \path[fill=black!70] (4,0) -- ++(1,0) -- ++(-1,1) -- cycle;
  \DrawDelta{5}{draw=none}{draw=none}
  \draw[line width=1.6pt] (2,0) -- (2,3);
  \fill (5,0) circle (2pt) node[below right] {$P_1$};
  \fill (0,0) circle (2pt) node[below left] {$P_0$};
  \fill (0,5) circle (2pt) node[above] {$P_2$};
  \node at (1.3,-0.5) {$X_1$};
  \node at (4.3,-0.5) {$X_4$};
  \node[right] at (2.2,2.9) {line $x=2$};
\end{tikzpicture}
\caption{The front half of the ring $X_1$ in $\T(5)$ (light gray) and of
the corner ring $X_4$ (dark gray). The proof of
Lemma~\ref{lem:no-mono} moves from $X_k$ to $X_{k+1}$ across the line
$x=k+1$ (thick), until it reaches the corner $P_1$.}
\label{fig:rings}
\end{figure}
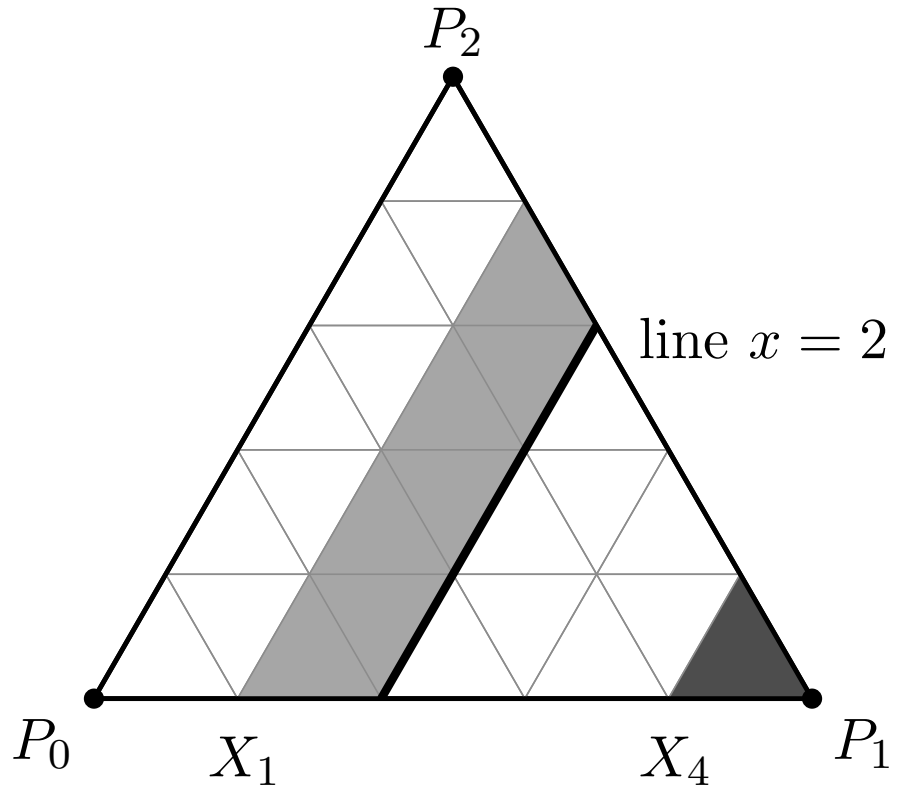

The ring $X_k$ closes up because $\Up^+(k,n-k-1)$ and $\Up^-(k,n-k-1)$
share their edge on the side $z=0$, and $\Up^-(k,0)$ and $\Up^+(k,0)$
share their edge on the side $y=0$. The ring $X_k$ has $2(2(n-k)-1)$
faces and the same number of crossing edges. Going once around $X_k$,
we pass from each face to the next one through exactly one crossing
edge, and we use every crossing edge exactly once. Inside a copy of
$\Delta_n$ this is the common edge of two consecutive triangles of the
strip; at the two ends it is the edge on the side $z=0$ or on the side
$y=0$. The last ring $X_{n-1}$ consists of the two faces
$\Up^\pm(n-1,0)$ at the corner $P_1=(n,0)$. These two faces share two
edges, and these are the two crossing edges of $X_{n-1}$. See
Figure~\ref{fig:rings}.

We define the $y$-rings $Y_k$ and the $z$-rings $Z_k$ in the same way,
using the strips $\{k\le y\le k+1\}$ and $\{k\le z\le k+1\}$. The
rotation of $\Delta_n$ by $2\pi/3$ about its center maps the lattice to
itself and permutes $x,y,z$ cyclically. Applied to both copies at
once, it induces a symmetry of $S_n$
that permutes the three families of rings. So anything we prove for
the $x$-rings also holds for the $y$-rings and the $z$-rings.

\begin{lemma}\label{lem:double-count}
Every edge of $S_n$ is a crossing edge of exactly two rings, and these
two rings belong to different families.
\end{lemma}

\begin{proof}
Let $e$ be an edge. Exactly one of $x,y,z$ is constant on $e$. Say $x$
is constant. Then $y$ takes the values $m$ and $m+1$ at the two ends
of $e$, for some integer $m$, so
$e$ lies in the strip $\{m\le y\le m+1\}$ and is a crossing edge of
$Y_m$. It is a crossing edge of no other $y$-ring. The same holds for
$z$. Finally, $e$ is a crossing edge of no $x$-ring, since $x$ does not
change along~$e$. If $e$ lies inside one copy of $\Delta_n$, it is a
crossing edge only for the rings that pass through that copy, and only
once. If $e$ lies on $\partial\Delta_n$, it is the edge where a ring
passes from the front to the back. In both cases $e$ is counted once in
each of the two rings.
\end{proof}

\subsection{No ring has a single color}

\begin{lemma}\label{lem:no-mono}
Let $\C$ be a good coloring of $\T(n)$. Then no ring is monochromatic.
\end{lemma}

\begin{proof}
By symmetry it is enough to treat the $x$-rings. By swapping colors, it
is enough to show that no $x$-ring is entirely black.

The ring $X_{n-1}$ consists of the two faces at the corner $P_1$, which
has degree~$2$. By Lemma~\ref{lem:local}(b) these faces have different
colors. So $X_{n-1}$ is not monochromatic.

Now let $0\le k\le n-2$. We claim that if $X_k$ is all black, then
$X_{k+1}$ is all black. Applying the claim $n-1-k$ times gives that
$X_{n-1}$ is all black, a contradiction. This proves the lemma.

To prove the claim, we look at the vertices on the line $x=k+1$. These
are the points $(k+1,j)$ with $0\le j\le n-k-1$. All of them have
degree~$6$, since none of them is a corner of $\Delta_n$. Each of them lies on faces of $X_k$
and of $X_{k+1}$ only.

\emph{The end point $(k+1,0)$.} It lies on $\partial\Delta_n$. On each
copy, it lies on the faces $\Up(k,0)$ and $\Dn(k,0)$ of the strip of
$X_k$, and on the face $\Up(k+1,0)$ of the strip of $X_{k+1}$. So four
of its six faces lie in $X_k$ and are black. By
Lemma~\ref{lem:local}(a), $b\in\{0,3,6\}$, so $b=6$. Hence
$\Up^+(k+1,0)$ and $\Up^-(k+1,0)$ are black.

\emph{Interior points.} Let $1\le j\le n-k-2$, and let $v$ be the point
$(k+1,j)$ on the front copy. Its six faces are
\begin{gather*}
  \Dn^+(k,j-1),\quad \Up^+(k,j),\quad \Dn^+(k,j)\qquad\text{in } X_k,\\
  \Up^+(k+1,j-1),\quad \Dn^+(k+1,j-1),\quad \Up^+(k+1,j)\qquad\text{in } X_{k+1}.
\end{gather*}
The first three are black. So $b(v)\ge3$, and Lemma~\ref{lem:local}(a)
gives $b(v)\in\{3,6\}$. Hence the last three faces all have the same
color.

Now we walk along the line. The face $\Up^+(k+1,0)$ is black. At the
point $(k+1,1)$ this forces $\Dn^+(k+1,0)$ and $\Up^+(k+1,1)$ to be
black. At $(k+1,2)$ this forces $\Dn^+(k+1,1)$ and $\Up^+(k+1,2)$ to be
black, and so on up to $j=n-k-2$. This covers every front face of
$X_{k+1}$. The same argument works on the back copy. If $k=n-2$ there
are no interior points, and $X_{n-1}$ consists only of
$\Up^\pm(n-1,0)$, which we already know are black. In all cases
$X_{k+1}$ is all black, which proves the claim.
\end{proof}

\subsection{Proof of Theorem~\ref{thm:main}}

Let $\C$ be a good coloring of $\T(n)$, and let $R$ be any ring. Walk
once around $R$, passing through each of its crossing edges once. The
color changes exactly when we pass through a crossing edge that is a
fold. Since we come back to the starting face, the number of
folds among the crossing edges of $R$ is even. By
Lemma~\ref{lem:no-mono} it is not zero. So it is at least~$2$.

There are $3n$ rings, namely $X_k$, $Y_k$ and $Z_k$ for $0\le k\le
n-1$. Let $N$ be the number of pairs $(R,e)$ where $R$ is a ring and $e$
is a fold that is a crossing edge of $R$. Counting by rings gives
$N\ge 2\cdot 3n=6n$. By Lemma~\ref{lem:double-count}, counting by folds
gives $N=2\fc(\C)$. Hence $\fc(\C)\ge3n$, and
\[
  \eta(\C)=\frac{\fc(\C)^2}{2n^2}\ \ge\ \frac{9n^2}{2n^2}=\frac92.
\]
Proposition~\ref{prop:B} shows that $\mathcal B_n$ attains this bound.
\qed

\begin{remark}
We also checked Theorem~\ref{thm:main} by computer for $1\le n\le8$.
We wrote the problem of minimizing the fold count as an integer
program, with one binary variable for each face and the congruence
of Definition~\ref{def:good} at each vertex. An exact solver returned
the minimum $3n$ in every case.
\end{remark}

\section{Subdivision}\label{sec:scaling}

\subsection{The subdivided triangulation}

Let $\beta\ne0$ be an Eisenstein integer and let $n\ge1$. Let $\mu\colon\Cx\to\Cx$ be the map $\mu(w)=nw$.
Then $\mu$ maps $\T$ onto the triangulation $n\T$ by triangles of side
$n$, and it maps $\beta\E$ onto $n\beta\E$. If $\rho$ is the rotation of
order~$3$ about a point $p$, then $\mu\rho\mu^{-1}$ is the rotation of
order~$3$ about $np$. So $\mu G_\beta\mu^{-1}=G_{n\beta}$, and $\mu$
induces an isomorphism
\[
  \T(\beta)=\T/G_\beta\ \longrightarrow\ (n\T)/G_{n\beta},
\]
which multiplies lengths by $n$.

Every triangle of $n\T$ is the union of $n^2$ triangles of $\T$. The
group $G_{n\beta}$ preserves both $\T$ and $n\T$, because it is generated
by rotations of order~$3$ about points of $n\beta\E\subset n\E\subset\E$,
and each of these maps both triangulations to themselves. So
$\T(n\beta)=\T/G_{n\beta}$ is obtained from $(n\T)/G_{n\beta}\cong\T(\beta)$
by cutting every triangle into $n^2$ triangles in the standard way. We
call the triangles of $\T(\beta)$ \emph{coarse} and the triangles of
$\T(n\beta)$ \emph{fine}. Each coarse edge is cut into $n$ fine edges.
The three vertices of degree~$2$ of $\T(n\beta)$ are images of points of
$n\beta\E$, so they are coarse vertices. All other points of
$\T(n\beta)$ have a flat neighborhood that looks like a piece of the
plane triangulation~$\T$.

For $\beta=1$ this is exactly the picture of
Section~\ref{subsec:doubled}. The coarse triangulation is $\T(1)$, with
two triangles, and $\T(n)$ is obtained by cutting each of them into $n^2$
unit triangles.

\subsection{Proof of Theorem~\ref{thm:scaling}}

Let $\C$ be a good coloring of $\T(\beta)$. Define $\C^{(n)}$ by giving
each fine face the color of the coarse face that contains it. We check
the condition of Definition~\ref{def:good} at every fine vertex~$v$.
There are three cases.

\emph{Case 1: $v$ lies in the interior of a coarse face.} Then $v$ has
degree~$6$, and all six fine faces at $v$ lie in the same coarse face.
So $b(v)\in\{0,6\}$, and the condition holds at $v$ by Lemma~\ref{lem:local}(a).

\emph{Case 2: $v$ lies in the interior of a coarse edge.} Then $v$ has
degree~$6$. Three fine faces at $v$ lie on one side of the coarse edge
and three lie on the other side. So $b(v)\in\{0,3,6\}$, and the condition holds at $v$.

\emph{Case 3: $v$ is a coarse vertex.} Each coarse face corner at $v$
contains exactly one fine face corner at $v$, and the two have the same
color. So $v$ has the same degree and the same numbers $b(v)$ and $w(v)$
in $\C^{(n)}$ as in $\C$. Since $\C$ is good, the condition holds at $v$.

So $\C^{(n)}$ is a good coloring of $\T(n\beta)$.

Now we count folds. If a fine edge lies in the interior of a coarse
face, its two fine faces lie in that coarse face and have the same
color. Otherwise the fine edge lies on a coarse edge~$e$, and its two
fine faces have the colors of the two coarse faces at~$e$. So the fine
edge is a fold if and only if $e$ is a fold of $\C$. Every coarse edge
contains exactly $n$ fine edges, so $\fc(\C^{(n)})=n\,\fc(\C)$. The
number of faces is multiplied by $n^2$. Therefore
\[
  \eta\bigl(\C^{(n)}\bigr)=\frac{n^2\fc(\C)^2}{n^2F}=\eta(\C). \qed
\]

\subsection{Proof of Corollary~\ref{cor:reduction}}

\begin{enumerate}[label=\textup{(\roman*)}]
\item This follows from Theorem~\ref{thm:main}: for every $n$ the coloring
$\mathcal B_n$ of $\T(n)$ has $\eta=9/2<\Omega$.

\item One direction is clear. For the other, assume that every primitive
$\T(\beta_0)$ has a good coloring with $\eta<\Omega$. Let
$\beta=a+b\alpha\ne0$ and let $n=\gcd(a,b)\ge1$ (with $\gcd(a,0)=|a|$). Then
${\beta_0=\beta/n}$ is a primitive Eisenstein integer and $\beta=n\beta_0$.
Take a good coloring $\C$ of $\T(\beta_0)$ with $\eta(\C)<\Omega$. By
Theorem~\ref{thm:scaling}, $\C^{(n)}$ is a good coloring of $\T(\beta)$
with ${\eta(\C^{(n)})=\eta(\C)<\Omega}$. \qed
\end{enumerate}

\section{Concluding remarks}\label{sec:remarks}

For the family $\{\T(n)\}$ we now know the exact minimum of $\eta$. The
lower bound $\eta\ge3$ of Theorem~\ref{thm:sch-lower} comes from
an isoperimetric inequality for each black or white region. Our bound
uses the closed rings of $\T(n)$ instead. It is sharp, but it uses the
special shape of $\T(n)$. For most other $\beta$ the edges of $\T$ are not parallel to the sides of the
doubled triangle, and the rings do not
close up in the same way.

By Corollary~\ref{cor:reduction}, Conjecture~\ref{conj:omega} is now a
question about primitive triangulations only. We end with two
questions.

\begin{question}
For each primitive $\beta$, let $\eta_{\min}(\beta)$ be the smallest
value of $\eta$ over all good colorings of $\T(\beta)$. Is
$\sup_\beta\eta_{\min}(\beta)$ finite?
\end{question}

\begin{question}
Up to the symmetries of $\T(n)$ and swapping colors, is $\mathcal B_n$
the only good coloring of $\T(n)$ with fold count $3n$?
\end{question}

\section*{Declarations}
\noindent\textbf{Conflict of interest.} The authors declare no conflict
of interest.\\
\textbf{Data availability.} No data were used. The computer code for
the check in Section~\ref{sec:lower} is available from the authors on
request.

\section*{Declaration of Generative AI and AI-Assisted Technologies in
the Manuscript Preparation Process}
During the preparation of this work, the authors used Claude AI
(Anthropic) to assist with organizing the presentation and improving the
language and readability of the manuscript. After using this tool, the
authors reviewed and edited the content as needed and take full
responsibility for the content of the published article.


\end{document}